\documentclass[12pt,a4paper]{article}
\usepackage{latexsym}
\usepackage{amssymb}
\begin{document}
\title{On quartics with three-divisible sets of cusps}

\author{S{\l}awomir Rams}
\date{{}}

\maketitle

\begin{abstract}
We study the geometry and codes of quartic surfaces with many cusps.
We apply Gr\"obner bases to find examples of various configurations of cusps on quartics.
\end{abstract}

\newcommand{\iiii}{\mbox{{\tiny i}}}
\newcommand{\iii}{\mbox{i}}
\newcommand{\epsi}{\varepsilon}
\newcommand{\q}{\mbox{q}}
\newcommand{\V}{\mbox{V}}
\newcommand{\X}{\mbox{X}}
\newcommand{\C}{\mbox{C}}
\newcommand{\Pic}{\mbox{Pic}}
\newcommand{\Fa}{\mbox{F}}
\newcommand{\Ha}{\mbox{H}}
\newcommand{\sing}{\mbox{sing}}
\newcommand{\supp}{\mbox{supp}}
\newcommand{\spann}{\mbox{span}}
\newcommand{\reg}{\mbox{reg}}
\newcommand{\mult}{\mbox{mult}}
\newcommand{\ov}{\overline}
\newcommand{\Aid}{A_{\mbox{id}}}
\newcommand{\Bid}{B_{\mbox{id}}}
\newcommand{\id}{\mbox{id}}
\newcommand{\ii}{\mbox{i}}
\newcommand{\pa}{\mbox{p}_{\mbox{{\footnotesize a}}}}
\newcommand{\pg}{\mbox{p}_{\mbox{{\footnotesize g}}}}
\newcommand{\ho}{\mbox{h}^{0}}
\newcommand{\hj}{\mbox{h}^{1}}
\newcommand{\hd}{\mbox{h}^{2}}
\newcommand{\calm}{\mbox{{\bf m}}}
\newtheorem{clm}{Claim}{\bf}{\it}
\newtheorem{case}{Case}{\bf}{\it}
\newtheorem{definition}{Definition}[section]
\newtheorem{lemma}{Lemma}[section]
\newtheorem{proposition}{Proposition}[section]
\newtheorem{cor}{Corollary}[section]
\newtheorem{thm}{Theorem}[section]
\newtheorem{example}{Example}{\bf}
\newcommand{\proofend}{\hspace*{\fill}$\Box$}
\newcommand{\proof}{\noindent{\it Proof}. }
\def\m{\mbox{{\bf m}}}
\def\kk{\mbox{{\takie K}}}
\def\CC{\mbox{{\takie C}}}
\def\NN{\mbox{{\takie N}}}
\def\PP{\mbox{{\takie P}}}
\def\QQ{\mbox{{\takie Q}}}
\def\ZZ{\mbox{{\takie Z}}}

\newcommand{\yczt}{Y_4}
\def\LL{{\cal L}}
\newcommand{\F}{{\mathbb F}_{\mbox{\tiny 3}}}
\newcommand{\KY}{{\cal K}_{\mbox{\footnotesize Y}}}
\newcommand{\KX}{{\cal K}_{\mbox{\footnotesize X}}}
\newcommand{\KW}{{\cal K}_{\mbox{\footnotesize W}}}
\newfont{\takie}{msbm10 scaled\magstep1}
\newfont{\taki}{msbm10}

\newcommand{\grad}{\mathop{\mathrm{grad}}\nolimits}
\newcommand{\ddj}{ \mbox{ $\frac{\partial Q}{\partial x_1}$ } }
\newcommand{\ddd}{ \mbox{ $\frac{\partial Q}{\partial x_2}$ } }
\newcommand{\ddt}{ \mbox{ $\frac{\partial Q}{\partial x_3}$ } }

\newcommand{\resyd}{X_d}
\newcommand{\sj}{S'}
\newcommand{\sd}{S''}
\newcommand{\Qd}{S}
\newcommand{\Ri}{R}
\newcommand{\Ei}{E'_{i}}
\newcommand{\Fi}{E''_{i}}
\newcommand{\Eo}{E'_{1}}
\newcommand{\Fo}{E''_{1}}
\newcommand{\LLp}{{\cal L}'}
\newcommand{\rr}{n}
\newcommand{\lr}{q}

\newcommand{\resyf}{X_4}
\def\FF{{\cal F}}


\section{Introduction}
The main aim of this note is to study 
the geometry and codes of quartics $Y_4 \subset  \PP_{3}( \CC ) $ 
with many cusps.  \\
 A {\em cusp}
(=singularity $A_2$) on $Y_4$ is a singularity near which
the surface is given in local (analytic) coordinates $x,y$ and $z$,
centered at the singularity, 
by an equation 
$$xy-z^3=0.$$
Let  $\pi:X_4 \to Y_4$ be the minimal desingularization
introducing two $(-2)$-curves $E_{i}', E_{i}''$ over each cusp $P_i$.
If 
there is a way to label these curves such that the divisor class of
$$\sum_{i =1}^n (E_{i}'+ 2 \, E_{i}'')$$
is divisible by three, 
then
the set $P_1,...,P_n$ is called three-divisible (\cite{bar1}, \cite{slt}).  
Equivalently: there exists a cyclic global
triple cover of $Y_4$ branched precisely over these cusps. In particular, every three-divisible set 
defines a vector (word) in the so-called code of the surface $Y_4$, see Sect.~\ref{codes}. 

It is well-known that a quartic surface  $Y_4 \subset  \PP_{3}( \CC ) $ has  at most eight cusps and 
a three-divisible set on $\yczt$ consists of six cusps (\cite{bar1}).
In \cite{slt} S.-L. Tan proves that
every quartic with eight
cusps contains a three-divisible set, but \linebreak

\vspace*{-0.5ex}

\footnoterule

\noindent
{\footnotesize
2000 {\it Mathematics Subject Classification.} 14J25, 14J17.

\vspace*{-0.2ex}
\noindent
{\it Key words.} quartic, A$_2$-singularity, ternary code.

\vspace*{-0.2ex}
\noindent
Research  partially supported by the Schwerpunktprogramm 
"Global methods in complex geometry" of the Deutsche Forschungsgemeinschaft. }

\pagebreak

\noindent
his 
method 
does not show which cusps on the surfaces defined in \cite{bar2} form such a set.
A general construction of quartics with three-divisible cusps is given in \cite{barr}.  
Here we use the results from  \cite{barr} to 
give a complete picture of the geometry of 
quartics with three-divisible sets:

In general the six cusps lie on a twisted cubic $C_3$, which however may degenerate 
to three concurrent lines in special cases. The surface $\yczt$
is given by the determinantal equation 
$$
 \det  \left [  \matrix{ S      &   Q_{12}     \cr 
                                    Q_{21} &   Q_{22} - S \cr } 
\right ]
\, = \, 0 , 
$$
where the quadrics $Q_{i,j}$ generate the ideal of $C_3$. 
 In particular, the locus of quartics with three-divisible sets 
is irreducible, but its germ   at each point corresponding to  a quartic with eight cusps
is reducible. \\
Every three-divisible set endows $\yczt$ with two elliptic fibrations. 
Their fibers are cut out by the entries of the above-given matrix,
and 
 multiplication of the matrix by a constant matrix corresponds to 
a choice of another pair of fibers. \\ 
We compute the code of a quartic with eight cusps and 
find all three-divisible sets on the quartics defined in \cite{bar2}. 
Finally, we use Gr\"obner bases to
give examples of both configurations of cusps.

In order to  render our exposition self-contained we recall some facts from \cite{barr}.
For various applications of three-divisible sets we refer the reader to \cite{slt}.   
We work over the field
of complex numbers.   


\section{Three-divisible sets and contact surfaces} \label{codes}

\noindent
Let 
$Y_4 \subset \PP_{3}( \CC ) $ be a quartic with finitely many singularities all of which are cusps
(i.e.~$\mbox{A}_2$-type double points).
Let  $\pi : X_4 \rightarrow \yczt$  be the  minimal resolution, 
where   $\sing(Y_4) = \{ P_{1}, \ldots,  P_{\lr} \}$ and
$E'_{i} + E''_{i} :=  \pi^{-1}(P_{i})$.

\vskip3mm
\noindent
{\bf Definition} (cf. \cite{bar1})  
{\sl If one can order the $(-2)$-curves $E'_{i},  E''_{i}$ in such a way that 
there exists a  divisor $\LL'$
which satisfies
$$
\sum_{1}^{6} (E'_{i} + 2  E''_{i}) = 3 \LL' ,
$$
then  $\{ P_{1}, \ldots,  P_{6} \}$ is called a $3$-divisible set.}

\vskip3mm
\noindent
Observe that 
if the divisor  $\sum_{1}^{6} (E'_{i} + 2  E''_{i})$ is $3$-divisible, then we can find an $\LL''$ such that 
$$
\sum_{1}^{6} (2 E'_{i} +   E''_{i}) = 3 \LL'' \, .
$$ 
Moreover, by \cite{br},  the divisor classes  $\LL'$, $\LL''$ are unique.

\noindent
In the sequel we assume that  $\yczt$  contains a $3$-divisible set $P_1, \ldots, P_6$,
and $\LL'$, $\LL''$  are  the above-defined divisors.  

We have the following lemma, which is proved in \cite{barr} under the 
assumption that the set $\sing( Y_4 )$ is $3$-divisible.
\begin{lemma} \label{nonplanar}
{\rm (cf. \cite[Lemma~3.2]{barr})} 
If $D \in |\pi^{*}{\cal O}_{Y_4}(1) - \LLp|$,
then $ \pi_{*}(D)$ is not a hyperplane section of  $Y_4$.
\end{lemma}
\noindent
\proof
Suppose to the contrary  and put $C :=  \pi_{*}(D)$.
Then we have 
\begin{eqnarray*}
 \pi^{*}{\cal O}_{Y_4}(1) - \LLp  & = &  
\tilde{C} + \sum_{1}^{\lr} ( \alpha'_i \,  \Ei +  \alpha''_i \, \Fi ) \, , \\
 \pi^{*}{\cal O}_{Y_4}(1) & = &  \tilde{C}  +
 \sum_{1}^{\lr} ( \beta'_i \, \Ei + \beta''_i \, \Fi ) \, ,
\end{eqnarray*}
where  $\alpha'_i$, $\alpha''_i$, $\beta'_i$,  $\beta''_i$ are 
non-negative integers 
and  $\tilde{C}$ is the  proper transform of the divisor $C$. Thus
$$
 \sum_{1}^{\lr} ( \alpha'_i \, \Ei +  \alpha''_i \, \Fi  ) +  \LLp =
 \sum_{1}^{\lr} ( \beta'_i \, \Ei +  \beta''_i \, \Fi )   \, . 
$$
We intersect that divisor with $(-2)$-curves $\Eo$, $\Fo$ to obtain the equalities 
$$ 
-2  \alpha'_1 + \alpha''_1 =  -2 \beta'_1 +  \beta''_1 \mbox{ and } 
\alpha'_1 - 2 \alpha''_1 - 1 =  \beta'_1 - 2  \beta''_1 \, . 
$$
Multiply the first equality by $2$ and add to the other to get 
contradiction with  the fact that  $\alpha'_1$,  $\beta'_1$ are integers. \proofend

We have $( \pi^{*}{\cal O}_{Y_4}(1) - \LL')^2 =    0$, so 
 Riemann-Roch yields :  
$$
\chi(  \pi^{*}{\cal O}_{Y_4}(1) - \LL')  = 2.
$$
Since the canonical system of $\resyf$ is trivial,
we can find divisors 
$$
D'_4 \in | \pi^{*}{\cal O}_{Y_4}(1) - \LL'| \mbox{ and }  D''_4 \in | \pi^{*}{\cal O}_{Y_4}(1) - \LL'')|.
$$
Let $C'_4 :=  \pi_{*}(D'_4)$ and let  $C''_4 :=  \pi_{*}(D''_4)$.
We have 
$$
3 \, C'_4 \sim  \pi^{*} {\cal O}_{Y_4}(3) - \sum_{i = 1}^{6} (E'_i + 2 E''_i),   
$$
so we can find a  cubic $\sj$  such that
$ \sj.Y_4 = 3 C'_4$.
Similarly, there is a cubic $\sd$ that satisfies   
$ \sd.Y_4 = 3 C''_4$.
Moreover, we can find a quadric  $\Qd$ such that  $\Qd.\yczt =  C'_4 + C''_4$. Indeed, we have 
$$
\pi_{*}{\cal O}_{X_4}((\pi^{*} {\cal O}_{Y_4}(1) - \LL') \otimes (\pi^{*} {\cal O}_{Y_4}(1) - \LL''))
 \subset  {\cal O}_{Y_4}(2) \otimes {\calm}_{P_{1}}  \otimes \ldots \otimes {\calm}_{P_{6}}.
$$ 
Having multiplied $\sj, \sd, \Qd$ by appriopriate constants, we get that
$(\sj \cdot  \sd)$ equals $(\Qd^{3})$ on $Y_4$, so     
the polynomial $(\sj \cdot \sd -  \Qd^3)$ vanishes on $\yczt$. 
One can prove that it does not vanish identically on  $\PP_{3}( \CC )$, see \cite[Thm~3.1]{barr}. \\
Let $R$ be the quadric  residual to $\yczt$ in the sextic
$(\sj \cdot \sd -  \Qd^3)$. Then   
$\yczt$ is given by the equation
$$
(\sj \cdot \sd -  \Qd^{3})/ R  = 0 \,  
$$  
and the cusps  $P_{1}, \ldots,  P_{6}$ belong to the surfaces $\sj, \sd, \Qd$. \\
In this section $\sj$,  $\sd$, $\Qd$ and  $R = R( \sj, \sd)$ denote the 
surfaces  constructed in this way for the quartic $\yczt$.


\begin{lemma} \label{lemmaongeneralfibres}
The system $|\pi^{*}{\cal O}_{Y_4}(1) - \LL'|$ has no base curve. Therefore,
the cubic  $\sj$ induced by push-forward of a general element of this system touches the surface
$\yczt$ along a smooth elliptic quartic with multiplicity three.
\end{lemma}
\proof 
Let $\sing( Y_4 ) = \{ P_1, \ldots, P_{\lr} \}$.
We 
put $\LL  :=  \pi^{*}{\cal O}_{Y_4}(1) - \LL'$.
Let $\FF$ be the free part of $\LL$ and let $C$ be a component of 
a general element of $|\FF|$. Then  $C^2 \geq 0$ (see \cite[p.~536]{gh})  and  $\deg( \pi(C) ) \leq 4$.
Since $\pa(C) = 1 + \frac{1}{2} \cdot C^2$, we have
$$
\pa( \pi(C))  \geq 1 \, .
$$
Hence $\pi(C)$ is either a quartic or a planar cubic. Thus if $L$ is a component of the base locus then either 
it is $\pi$-exceptional or $\pi(L)$ is a line. Moreover, at most one fixed component of $|\LL|$ is not $\pi$-exceptional.  

\noindent
Case 1: Let us assume that all fixed components are   $\pi$-exceptional, i.e.
$$
\LL = \FF +  \sum_{1}^{\lr} (  \alpha'_i \, E'_i + \alpha''_i \, E''_i) \, .
$$
Then  $\FF. C \geq 0$ for every irreducible $C$. From the equalities
$$
0 = \LL^2 = \LL. ( \FF +  \sum_{1}^{\lr} (\alpha'_i \, E'_i +  \alpha''_i \, E''_i)) = 
\LL. \FF + \sum_{1}^{6}  \alpha''_i \, ,
$$  
we obtain  $\LL. \FF =  0$. Moreover, we have 
$\LL.  E'_i =  0$ for $i \leq \lr$ and  $\LL.  E''_i = 1$  (resp.  $\LL.  E''_i = 0$) for $i \leq 6$ (resp.  $i > 6$),
which  implies
$$
0 = (\FF +  \sum_{1}^{\lr} (  \alpha'_i \, E'_i + \alpha''_i \, E''_i)). \FF  \geq  \sum_{1}^{\lr} ((\alpha'_i)^2 +  (\alpha'_i -  \alpha''_i)^2 + (\alpha''_i)^2)  \,  .
$$

\noindent
Case 2: Let $\LL = \FF +  L + \sum_{1}^{\lr}( \alpha'_i \, E'_i + \alpha''_i \, E''_i)$, where 
$\pi(L)$ is a line.
Then 
$$
\LL . L = ( \pi^{*}{\cal O}_{Y_4}(1) - \frac{1}{3} \sum_{1}^{6} (E'_{i} + 2 E''_{i})).L =
 4 -   \frac{1}{3} \sum_{1}^{6} (E'_{i}.L + 2 E''_{i}.L) \, .
$$
The latter yields  $\LL . L > 0$.  Indeed,
since $E'_{i}.L \leq 1$ and $E''_{i}.L \leq 1$, the inequality $\LL . L \leq 0$ implies that at least four cusps lie on $\pi(L)$, so 
$\yczt$ is singular along $\pi(L)$. Contradiction. We get
$$
\LL^2 = ( \FF + L +  \sum_{1}^{q} (\alpha'_i E'_i +  \alpha''_i E''_i)). \LL = 
\FF . \LL + L. \LL + \sum_{1}^{6} \alpha''_i > 0 \, ,
$$
which contradicts $\LL^2 = 0$. \\
Thus general element of $|\LL|$ is irreducible and, by 
$\LL^2 = 0$,  it is smooth. 

Finally,  by Lemma~\ref{nonplanar} the curve  $\pi(C)$ is not planar, so it is a quartic. 
In particular $\pa( \pi(C)) = 1$ (see \cite{har}).
If $\pi(C)$
were singular, 
then its proper transform under the blow-up of cusps of $\yczt$ 
would be either singular or smooth rational. This is in conflict with 
the fact that  $C$ is smooth elliptic. Thus general
$\sj$ meets $\yczt$ along a smooth elliptic quartic. 
\proofend


\section{The equation of $\yczt$} \label{theequation}

\noindent
Here  we find another equation of  $\yczt$. The proof of the main proposition is preceded
by two lemmas. 
We maintain the notation and the assumptions of Sect.~\ref{codes}. Moreover,
in  Lemmas~\ref{lemmaonrestrictionstoresidual},~\ref{residualisreduced} we assume that 
\begin{itemize}
\item the quadric $\Qd$ is smooth,
\item the quartic curves $C'_4 := \sj \cap \yczt$, $C''_4 := \sd \cap \yczt$
 are irreducible.
\end{itemize}
The latter implies that the cubics $\sj$, $\sd$ are irreducible. 

\begin{lemma} \label{lemmaonrestrictionstoresidual} 
There exist (possibly non-reduced) conics $C', C''$ such that  
$$
\sj. R = 3 \,  C', \, \, \, \, \, \,  \sd.R = 3 \, C'' \mbox{ and } \Qd \cap R = C' \cup  C'' \, .
$$  
\end{lemma}
\proof 
We prove that every component appears in the cycle
$\sj . R$ with a  $3$-divisible multiplicity. 
Let $S_6 := (\sj \cdot  \sd - \Qd^3)$. Then
$$
 \sj.S_6 =  \sj .  R + \sj . \yczt =  Z'_{1} + Z'_{2} +  3 C'_4,
$$
where  $Z'_{1}$ is the part of the cycle $ \sj .  R$ with no  components in 
$\sing( \sj)$ and $\mbox{supp}(Z'_{2}) \subset \sing( \sj)$. 
Since
$S_6 \equiv - \Qd^3$ on $\sj$,  every $C \not \subset \sing( \sj)$ appears in  
$S_6 . {\sj}$ with a $3$-divisible multiplicity. In particular,
 $\deg(Z'_{1})$  is divisible by $3$. 
If $Z'_{2} = 0$, then we are done. Otherwise,  
 $\supp(Z'_{2})$ is a line. Its multiplicity in the cycle $Z'_{2}$ is  $\deg(Z'_{2})$, so it is 
  $3$-divisible.

\noindent
Let $C' :=   \supp ( \sj . R )$. Since $\deg( C' ) \leq  2$, it suffices to prove that $C'$ is connected. 
This is obvious when $R$ is irreducible, 
so we can assume  that $R$ is a union of planes and $C'$ consists of two lines. 
Then $S$, $S'$ intersect along the irreducible quartic $C'_4$ and the lines, so the latter 
form a $(1,1)$-curve.  
The proof for $\sd$ and $C''$ is analogous. 
 
The quadrics
 $\Qd$, $R$ meet along the curves $C',C''$ because 
$\sj \cdot  \sd \equiv  \Qd^3$ on   $R$ ( 
if  $R$ is a double plane, then  consider  restriction  to 
this plane).
\proofend

\begin{lemma} \label{residualisreduced}
The quadric $R$ is reduced and 
 the conics 
$C'$, $C''$ are its hyperplane sections. Moreover, if
 $R = R_1 \cdot R_2$ consists of two planes, then
\begin{equation} \label{hyperplanesection}
C' = L.R ,
\end{equation}
where the forms $L$, $R_1$, $R_2$ are linearly independent.
\end{lemma}
\proof 
We  maintain the notation of Lemma~\ref{lemmaonrestrictionstoresidual} 
 and prove that $R$ is reduced. \\
Suppose that  $R = x_0^2$. The sextic $\sj \cdot \sd$ meets $R$ along two distinct lines, say $l', l''$,
because $S$ is smooth.
By Lemma~\ref{lemmaonrestrictionstoresidual} we have $\sj.R = 6 l'$ and   $\sd.R = 6 l''$, so
we can assume that 
\begin{center}
\begin{tabular}{lll}
 $\sj$ & = & $ x_1^3 + x_0 \cdot ( x_0 \cdot S'_1 +   S'_2)$  , \\
 $\sd$ & = & $ x_2^3 + x_0 \cdot ( x_0 \cdot S''_1 +   S''_2)$ , \\
 $\Qd$ & = & $ x_1 \cdot x_2 +  x_0 \cdot ( x_0 \cdot s_1  +   S_2)$ , 
\end{tabular}
\end{center}
where $S'_2, S''_2, S_2  \in \CC[x_1, x_2, x_3]$. Then 
$$
S_6 = 
 x_0  \cdot ( x_2^3 \cdot S'_2  + x_1^3 \cdot S''_2  -  3 \,  x_1^2 \cdot x_2^2 \cdot S_2) = 0 \, \, \, \, 
(\mbox{mod } x_0^2), 
$$
so
 $S'_2$ (resp. $S''_2$) is divisible by  $x_1^2$ (resp. $x_2^2$),
and $S_2 \in  \CC[x_1, x_2]$. 
The latter  is in conflict with the assumption that $\Qd$ is smooth. 

Suppose that  $R$ is irreducible. If $R$ is a cone, then  $C' \in {\cal O}_R(1)$.
Otherwise, $\sj .R = 3 C'$ is a $(3,3)$-curve, so
 $C'$ is a  hyperplane section of $R$. \\
Assume that $R = R_1 \cdot R_2$ consists of two planes. Then
$C' = l_1 + l_2$ is a sum of lines.
Suppose that the planes $R_1$,  $R_2$ meet along the line $l_1$. Then $l_1 = l_2$. Indeed, 
assume that $l_2 \subset R_1$ and observe that
if  $l_2 \neq l_1$,  then by Lemma~\ref{lemmaonrestrictionstoresidual} the cycle $\sj.R_1 -  3 l_2 - l_1$ 
is non-negative, so $\deg (\sj.R_1) \geq 4$. This is in conflict with the irreducibility of $\sj$. 
Thus $S'$, $R$ 
meet only along $l_1$. Since the conics $S.R_1$, $S.R_2$ contain $l_1$, the cubic $S''$ meets 
$R$ in two distinct lines $l''_1$, $l''_2$, none of which coincides with $l_1$. Thus 
$$
S.R = 2 l_1 + l''_1 + l''_2 \, ,
$$ 
so $2 l_1$ is  a $(1,1)$-curve. Contradiction. 
\proofend

\begin{proposition} \label{equationsofquartics} 
Let $\yczt$ be a quartic with a $3$-divisible set
$\{ P_1, \ldots, P_6 \}$ of cusps.
 Then there exist linear forms  $L',L'',F',F''$
 and   a quadric $R$  such that $\yczt$ is given by the equation: 
\begin{equation} \label{minus}
(\sj \cdot \sd - \Qd^3)/ R = 0,
\end{equation}
where  the forms  $L',  L''$ are linearly independent,
\begin{equation} \label{kubiki}
\sj := (L')^3 + F' \cdot R, \, \, \sd := (L'')^3 + F''\cdot R \mbox{ and }
\Qd :=  R + L' \cdot  L''. 
\end{equation} 
Moreover,  the inclusion
$\{P_1, \ldots, P_6 \} \subset  \sj \cap \sd \cap \Qd$ holds.
\end{proposition}
\proof
We claim  that $\yczt$ is given by the equation (\ref{minus})
with a smooth quadric $\Qd$. Indeed,
general cubic
$\sj$ meets $\yczt$ along a smooth irreducible quartic $C'_4$. 
Once we fix $\sj$ and vary $\sd$, the quadric $\Qd$ varies but always contains the curve  $C'_4$. 
There are only finitely many singular quadrics that contain $C'_4$, so  the quadric $\Qd$ is 
smooth for general $\sj, \sd$. \\
Lemma~\ref{residualisreduced} yields that  the residual quadric $R$ is reduced. Assume that    
 $R = R_1 \cdot R_2$ consists of two planes and choose a form $L$ which satisfies (\ref{hyperplanesection}). Then the function 
 $(\sj / L^3)$ is constant on  both planes $R_i$, so 
$$
\sj = \alpha' \cdot  L^3 +  R_1 \cdot  Q'_1 
=  \alpha'  \cdot L^3 +  R_2  \cdot Q'_2 = 
 \alpha' \cdot L^3 +  R_1 \cdot R_2 \cdot F',
$$   
where the last equality results from the fact that  $L, R_1, R_2$ are linearly independent.
Moreover, $\alpha'$ is non-zero because the cubic $\sj$ is irreducible, so we can assume  $\alpha' = 1$. \\
In the same way we find forms $L'', F''$ 
such that  
$\sd$ satisfies (\ref{kubiki}). \\
The function $( S / (L' L''))$ is constant on  $R_i$
 because  $S^3 \equiv   (L' L'')^3$ on that planes. Therefore,
$$
\Qd = \alpha \cdot  L' \cdot L'' + \beta \cdot R \mbox{ with }  \alpha^3 = 1.
$$
Multiply  $L'$, $F'$, $F''$ and $R$ by appropriate constants to obtain (\ref{kubiki}). \\
The proof of  (\ref{kubiki}) for an irreducible quadric $R$ follows the same lines , so we leave it to the reader.

The  forms  $L', L''$ are  linearly independent because the polynomial 
\begin{equation} \label{qij}
\sj - F' \cdot S = L' \cdot ((L')^2 -   F' \cdot L'')
\end{equation}
 belongs to the ideal of 
the non - planar irreducible curve $C'_{4}$.
\proofend

\section{Configuration of cusps}


\noindent
Here we study  the configuration of cusps in $3$-divisible sets.
We introduce the following notation:
\begin{center}
\begin{tabular}{lll}
 $Q_{1,2}(L', \ldots,  F'')$ & := & $L' \cdot F'' - (L'')^2$, \\
 $Q_{2,1}(L', \ldots,  F'')$ & := & $L''\cdot F' - (L')^2$,   \\
 $Q_{2,2}(L', \ldots, F'')$  & := & $F' \cdot F'' - L' \cdot L''$, \\
 $C_3(L', \ldots, F'')$      & := & $Q_{1,2} \cap Q_{2,1} \cap Q_{2,2}$ \, .
\end{tabular}
\end{center}
We omit the linear forms  and 
write $Q_{i,j}$, resp. $C_3$ when it is  not ambiguous. 
\begin{thm} \label{positionofcusps}
Every quartic $\yczt$ with a $3$-divisible set
 is given 
 by the equation:
\begin{equation} \label{determinantalequation}
 \det  \left [  \matrix{ S      &   Q_{12}     \cr 
                                    Q_{21} &   Q_{22} - S \cr } 
\right ]
\, = \, 0 , 
\end{equation}
where $S$ is  smooth and the $3$-divisible set
consists of the intersection points of the quadric $\Qd$ and the curve $C_3$. Moreover, one of the following holds:

\vskip2mm
\noindent
{\rm ({\bf I})}  The forms $L', L'', F', F''$ have no common zero, i.e. $C_3$ is a twisted cubic.

\vskip2mm
\noindent
{\rm ({\bf II})} The planes $L', L'', F', F''$ meet in one point $P$. Then $C_3$ consists of
 three 
lines that meet in $P$. Each of the lines contains precisely two cusps. The vertex $P$ does not lie on the 
quartic $\yczt$ and the lines are not coplanar.
\end{thm}
\proof
Let $P_1, \ldots, P_6$ be the $3$-divisible set on $\yczt$.
The equation (\ref{determinantalequation}) can be obtained from (\ref{minus}) by direct computation.  
Furthermore, since
$$
\{ P_1, \ldots, P_6 \} \subset  \sj \cap \sd \cap \Qd \, ,
$$
all the cusps $P_1, \ldots, P_6$ belong to 
the quadrics $Q_{1,2}, Q_{2,1}$ (see (\ref{qij})).
It remains to compare multiplicities in a cusp  $P_j$:
$$
2 = \mult_{P_j}( \yczt ) = \mult_{P_j}( \Qd \cdot (Q_{2,2} - \Qd )) = 1 +  \mult_{P_j}(Q_{2,2} - \Qd) \, , 
$$
so $P_j$ lies on $Q_{2,2}$, and  the cusps $P_1, \ldots, P_6$ belong to $\Qd \cap C_3$.

If  $L', L'', F', F''$ have no common zero, then $C_3$ is  twisted cubic and  we have
 $ \Qd.C_3 = \sum P_j$. This is the configuration (I). 

Assume that the planes  $L', \ldots, F''$ meet.
Suppose that  their intersection is a line $l$.
Then the quadrics  $Q_{i,j}$ are singular along that line, so they are sums of planes.
If they met  along a surface, then $\yczt$ would be singular along its intersection with 
the quadric $\Qd$, so it would be singular along a curve. Hence the quadrics  $Q_{i,j}$ meet precisely 
along the line $l$ and  the $3$-divisible cusps
lie on that line. This is impossible, so
 the forms  $L', \ldots, F''$ have  only one common zero, say $P$. \\
Since $Q_{i,j}$ are cones with vertex $P$, the curve $C_3$ is a cone,  i.e. a sum of lines 
that pass through $P$. 
The lines  are not coplanar by Lemma~\ref{lemmaongeneralfibres}, so there are  at least three of them. 
If $Q_{1,2}$, $Q_{2,1}$
have no common components, then their intersection is a quartic curve with the line $L' \cap L''$ as a component.
The latter does not lie on $Q_{2,2}$, hence
$C_3$ consists of three lines. Otherwise, the quadrics  $Q_{1,2}$, $Q_{2,1}$ meet along
 a line  and  a plane. The latter intersects $Q_{2,2}$ properly,  so 
$\deg( C_3 ) = 3$ again. \\
Finally, the quartic $\yczt$ meets $C_3$ in  points on  $\Qd$,  i.e. in the cusps. Since there are six of them, 
every component of $C_3$ intersects $\Qd$ in two cusps and      
 the vertex $P$ does not belong to  $\yczt$. 
\proofend 

Let ${\bf {\cal D}_4}$ denote the closure of the  family of quartics given by the equation 
$$
 \Qd \cdot ( Q_{2,2} - \Qd) - Q_{1,2} \cdot Q_{2,1} = 0 \, ,
$$
where $\Qd$ is a quadric,  $L', L'', F', F''$ are 
linear forms,
 $Q_{i,j} := Q_{i,j}(L', \ldots, F'')$. \\
 By  Thm~\ref{positionofcusps} every  quartic with a $3$-divisible set
  belongs to  ${\cal D}_4$. Thus \cite[Thm~2.1]{barr} 
implies that 
 projective quartics with $3$-divisible sets of cusps 
form a dense subset of ${\cal D}_4$. 
\begin{cor}
The locus of quartics with $3$-divisible sets is irreducible. 
\end{cor}

\newcommand{\ta}{\underline{a}}
\newcommand{\tA}{\mathfrak{a}}
\newcommand{\ttaj}{A_1}
\newcommand{\ttad}{A_2}
\vskip3mm
It is natural to ask how to pass from the equation (\ref{determinantalequation})
to (\ref{minus}), i.e. how to find the cubics that touch $\yczt$ with multiplicity $3$.
Here we answer this question in generic case, i.e. for type (I). The answer in the other case is similar.  \\
Let  $\yczt$ be a quartic carrying a $3$-divisible set 
of type (I) with  $L' := x_0, \ldots, F'' := x_3$ and let 
$\Phi ( t_0, t_1 ) := ( t_0^2 \, t_1 , t_0 \,  t_1^2 , t_0^3 , t_1^3 ).$ \\
For a non-degenerate $(2 \times 2)$-matrix $\ta := [ a_{i,j} ]_{i,j = 0,1}$
there exists a unique automorphism 
$\tA = ( L''( \ta ), L'( \ta ), F''( \ta ), F'( \ta ) )$
of $\CC^4$  such that $\tA \circ \Phi = \Phi \circ \ta$. Put
$$
\ttaj :=    \left  [  \matrix{  a_{0,1}     &   a_{0,0} \cr 
                      a_{1,1} & a_{1,0} \cr } \right ], \, \, \, \, \, \, \, \,  
\ttad :=    \left  [  \matrix{  a_{1,0}     &   a_{0,0} \cr 
                      a_{1,1} & a_{0,1} \cr } \right ] \,
$$
Then, by direct computation,  one can find a quadric $Q(\ta)$ such that 
$$
 \ttaj \, 
\left [  \matrix{ Q      &   Q_{12} \cr 
                      Q_{21} & Q_{22} - Q \cr }  \right ] \, \ttad =
 \det( \ta )^{(-2)} \left [  \matrix{ Q(\ta)      &   Q_{12}(\ta) \cr 
                      Q_{21}(\ta) & Q_{22}(\ta) - Q(\ta) \cr }  \right ] \, , 
$$
where $Q_{i,j}(\ta) := Q_{i,j}(L'(\ta), \ldots, F''(\ta))$. 
Now it suffices to use (\ref{kubiki}) to find the cubics. \\
The above equality  shows that
the intersection of a cone over $C_3$ (with vertex on $C_3$) with the quartic $\yczt$ consists of 
two quartic curves, each of which is image under $\pi$ of a curve in   $|\pi^{*}{\cal O}_{Y_4}(1) - \LL'|$, resp.
$|\pi^{*}{\cal O}_{Y_4}(1) - \LL''|$. The choice of the vertex of the cone corresponds to the choice of a fiber in the elliptic fibration 
given by $|\pi^{*}{\cal O}_{Y_4}(1) - \LL'|$.  
  
\section{Quartics with eight cusps} \label{eightonquartics}

Here we study quartics with eight cusps. In particular we find $3$-divisible sets on 
the quartics defined in \cite{bar2}. 
\begin{definition}
The code of the quartic $\yczt$ is the kernel of the map 
$$
\varphi:  \F ^{q} \ni  \sum_{1}^{q} \mu_{j} \, e_{j} \mapsto \sum_{1}^{q} \mu_{j} \,  (E'_{j} - E''_{j})  \, \in  
{\rm Pic}(X_4) \otimes \F .  
$$
\end{definition}
Observe that every vector (word) in  the linear code $\mbox{ker}(\varphi)$ corresponds to 
a~$3$-divisible set of cusps on $\yczt$.  

\noindent
For   $v \in \F^{q}$, one defines its weight 
as the number of its non-zero coordinates.
A ternary $[q, d, \{ r \}]$ {\sl code} is a  $d$-dimensional subspace $\F^{q}$ such that 
all its non-zero elements are of weight $r$.
We have the Griesmer bound \cite[Thm~(5.2.6)]{vlt}:
\begin{equation} \label{gb}
q \geq \sum_{i = 0}^{i = d-1} \lceil \mbox{{\small $\frac{r}{3^{i}}$}} \rceil \, 
\end{equation}

\begin{thm} \label{maximalnumberofcuspsonquartics} 
If $\yczt$ is a  quartic with eight  cusps, then 

\noindent
(i) the germ of ${\cal D}_{4}$ at $\yczt$ is reducible,

\noindent
(ii) the  code of  $\yczt$ is the  $[8, 2, \{ 6 \}] - \mbox{code}$ that  is spanned by the following words  
$$
\begin{tabular} {c c c c c c c c   }
 1 & 1  & 1 & 1  & 1 & 1 & 0  & 0   \\
 0 & 0  & 1 & 1  & -1  &  -1 & 1 & 1  
\end{tabular}
$$ 
\end{thm}
\proof
(ii) The above-defined code is the only ternary  $[8,2, \{ 6 \}]-\mbox{code}$.
Moreover, by (\ref{gb}), every  
 $[8, p,  6]$-code is at most $2$-dimensional, so it suffices to prove that 
$\mbox{ker}( \varphi )$ is $2$-dimensional.  
Suppose that $\dim \mbox{ker}( \varphi ) \leq 1$. Then $\dim( \mbox{Im}( \varphi)) = 15$ and 
 $\dim( \mbox{Im}( \varphi)^{\perp}) = 7$. 
On the other hand we have 
$$
E'_i \in  \mbox{Im}( \varphi)^{\perp} \mbox{ and } (E'_i)^2 = (-2) \, , 
$$ 
so  $\dim( \mbox{Im}( \varphi)^{\perp}) \geq 8$. Contradiction. \\
(i) Consider the projection from the set 
$$
\{ (x, Y) :   x \in \sing( Y ), Y \in {\cal D}_{4}  \}
$$
on the variety ${\cal D}_4$. By \cite{barr} its general fiber consists of six points, but the fiber over $\yczt$ has eight elements. 
This cannot happen when the germ of  ${\cal D}_{4}$ at $\yczt$ is irreducible.
\proofend

\vskip2mm
By  Thm~\ref{maximalnumberofcuspsonquartics} 
two 
 distinct $3$-divisible sets never have  five cusps in common.

\vskip2mm
\noindent
{\bf Example  \ref{eightonquartics}.1} {\rm (cf. \cite{bar2})} Let $k \neq 0$. The surface $S_k$ given by the polynomial 
\begin{eqnarray*}
\begin{array}{l}
(1 + k)^3 x_0^2 x_1^2 + 2 k (1 - k^2) x_0 x_1 x_2 x_3 - (1- k)^3 x_2^2 x_3^2 \hfill \\
+  (1- k)^2 (x_0 +  x_1 + x_2 +  x_3) [ (1 - k) x_2 x_3 (x_0 +  x_1) -  (1 + k) x_0 x_1 (x_2 +  x_3) ]
\end{array}
\end{eqnarray*}
has the obvious $\ZZ_2 \times \ZZ_2$-symmetry:
$$
\psi_1 : x_0  \leftrightarrow x_1, \, \, \,  \psi_2 : x_2  \leftrightarrow x_3.
$$
The set of its cusps splits into three orbits of the  $\ZZ_2 \times \ZZ_2$-action:
\begin{center}
\begin{tabular} {l l c  l l  c l }
$({\mathfrak p}_{1})$   & $P_1$ & $:=$ & $(1 : 0 : -1 : 0)$, & $ P_{2}$ & $:=$ & $(1 : 0 : 0 : -1)$, \\
        & $P_{3}$ & $:=$ & $(0 : 1 : -1 : 0)$, & $ P_{4}$ & $:=$ & $(0 : 1 : 0 : -1)$, \\
$({\mathfrak p}_{2})$ & $P_5$     & $:=$ & $(1 : 0 : 0 : 0)$,  & $ P_6    $ & $:=$ & $(0 :1 : 0 : 0)$,   \\
$({\mathfrak p}_{3})$ & $P_7$     & $:=$ & $(0 : 0 : 1 : 0)$,  & $ P_8    $ & $:=$ & $(0 :0 : 0 : 1)$.   
\end{tabular}
\end{center}  
Let ${\cal A}$ be a $3$-divisible set on $S_k$. Suppose that ${\cal A}$ contains precisely three elements of the orbit 
$({\mathfrak p}_{1})$  . W.l.o.g.  $P_{1}, P_{2}, P_{3}$ belong to the $3$-divisible set ${\cal A}$.
If either $P_5 \not \in {\cal A}$  or  $P_6 \not \in {\cal A}$, then $\psi_1( {\cal A} )$ and 
${\cal A}$ have five cusps in common, which is in 
conflict with  Thm~\ref{maximalnumberofcuspsonquartics}.
Similar reasoning applied to  $\psi_2$ rules out the remaining possibilities, i.e.  $P_7 \not \in {\cal A}$, 
$P_8 \not \in {\cal A}$. \\
Assume that  ${\cal A}$ contains the orbit $({\mathfrak p}_{1})$. As before   
Thm~\ref{maximalnumberofcuspsonquartics} implies that 
 ${\cal A}$  contains one of the orbits $({\mathfrak p}_{2})$, $({\mathfrak p}_{3})$. 
Finally, let  ${\cal A}$ contain the orbits  $({\mathfrak p}_{2} )$, $({\mathfrak p}_{3})$.
By Thm~\ref{positionofcusps} no five cusps in  ${\cal A}$ are coplanar, 
so the $3$-divisible set contains either the points $P_{1}, P_{4}$ or 
 $P_{2}, P_{3}$.
To sum up,  ${\cal A}$ is one of the following 
\begin{eqnarray*}
\begin{array} {c }
 P_1, \ldots, P_4, P_5, P_6 \\
 P_1, \ldots, P_4, P_7, P_8 \\
 P_5, \ldots, P_8, P_1, P_4 \\
 P_5, \ldots, P_8, P_2, P_3 
\end{array}
\end{eqnarray*}   
According to  Thm~\ref{maximalnumberofcuspsonquartics} there are four $3$-divisible sets on $S_k$, so all the 
above-given sets are  $3$-divisible.
They are of type (II) because each of them contains four coplanar points.

\section{Examples} \label{examples}

\noindent
In order to detect  a  $3$-divisible set of cusps we use the following lemma.

\begin{lemma} \label{genaraldescriptionofthreedivisibilty} 
Let $\yczt$ be a quartic carrying a set $P_1, \ldots, P_6$ of cusps such
that it has no other singularities.
Let $\sj$ be a cubic that touches $\yczt$ with multiplicity three, e.g. $\sj.\yczt = 3 C'_4$, and  $P_1, \ldots, P_6 \in \sj$. 
If the divisor $C'_4$ is smooth in the cusps, then they form a $3$-divisible set. 
\end{lemma}

\proof 
Let   $\alpha'_{i} \leq  \alpha''_{i}$ be the positive integers that satisfy  
$$
\pi^{*}\sj = 3 \tilde{C}'_4 + \sum_{1}^{6} ( \alpha'_{i} \,  E'_{i} +  \alpha''_{i}\,  E''_{i}).
$$
We claim that $ \alpha'_i = 1 \mbox{ and }  \alpha''_i = 2  \mbox{ for } i = 1, \ldots, 6$. Indeed, let 
$H \subset  \PP_{3}$ be a general hyperplane  through the point $P_i$. In particular, 
we assume that
$$
\pi^{*}H = \tilde{H} + E'_{i} +  E''_{i} \, .
$$
W.l.o.g the plane $H$ contains no components of the tangent cone of the support of the divisor  
$C'_4$ at $P_i$, which implies that the proper transforms
$\tilde{H}$, $\tilde{C}'_4$ meet in no points on the exceptional divisors $E'_i$, $E''_i$.  
Moreover we have $\tilde{H}.\tilde{C}'_4 = (\deg(C'_4) - \mult_{P_{i}}(C'_4)) = 3$. 
It results from  $(\pi^{*}S').E'_i =   (\pi^{*}S').E''_i = 0$ that 
$(3 \tilde{C'_4}).(E'_i + E''_i) =  \alpha'_i + \alpha''_i$.
Our claim follows from the equality
$$
12 = (\pi^{*}\sj).(\pi^{*}H) =
 (3 \,  \tilde{C'_4}).( \tilde{H} + E'_{i} +  E''_{i}) = 9 +  \alpha'_i + \alpha''_i \,.
$$
Finally, we have
$\sum_{1}^{6} (  E'_{i} + 2  E''_{i})   =
3 (\pi^{*} {\cal O}_{Y_4}(1) -  \tilde{C'_4})$.
\proofend

\vskip2mm
We give examples of both configurations of cusps (see Thm~\ref{positionofcusps}).


\vskip2mm
\noindent
{\bf Example \ref{examples}.1} {\rm (Configuration (I))} We put
\begin{center}
\begin{tabular}{lll }
$\Qd$ & := & $49 x_1^2 + x_2^2 - 36 x_3^2  - 14 x_0^2$, \\
$R$   & := & $\Qd - x_0 \cdot x_1$, \\
$\sj$ & := & $x_0^3 + x_2 \cdot R$, \\
$\sd$ & := & $x_1^3 + x_3 \cdot R$,
\end{tabular}
\end{center}
and define the quartic $\yczt$ by the equation (\ref{minus}).
We claim that $\sing( \yczt )$  consists of the six cusps 
$(\pm j : j : j^3 : \pm 1 ) \, , \mbox{ where } j = 1, 2, 3,$
that form a $3$-divisible set of type (I). \\
By direct computation   $\sj$, $\sd$, $\Qd$ meet transversally in the six points, none of which 
belongs to  $R$,   
so those points  are cusps on  $\yczt$. If we put 
$$
L' := x_0, \, L'' := x_1, \,  F' := x_2, \,  F'' := x_3,
$$
 then 
$\yczt$ is given by the equation (\ref{determinantalequation})
and the curve $C_3$ is  twisted cubic. 
A~Gr\"obner basis computation (see Remark~\ref{examples}.1) shows that the polynomials 
$ (Q_{i,j})^4, \Qd^4$ belong to the jacobian ideal of $\yczt$, so the 
 only singularities of  $\yczt$  are the six cusps.
It remains to  prove that they form a $3$-divisible set.

Observe that the curves   $C'_4 := \yczt \cap \sj$,  $C''_4 := \yczt \cap \sd$ 
have no common components. Indeed, we have
$$
S' \cdot S'' - S^3 = \yczt \cdot R,
$$
so, by multiplicity count,  a smooth point of $\yczt$ that 
belongs to both curves must lie on $R$. But, the surfaces  $R$, $\sj$, $\sd$ meet properly because 
the line $x_0 = x_1 = 0$ does not lie on the quadric $S$. 
Thus the cubic  $\sj$ touches $\yczt$ with multiplicity $3$ along  $C'_4$.
It suffices to prove that $C'_4$ is smooth in the cusps (see Lemma~\ref{genaraldescriptionofthreedivisibilty}).

Let $H$ be a general hyperplane 
through a fixed cusp $P_i$ 
and let $\mbox{i}( \,  \cdot \,  )$ stand for the intersection multiplicity. Then
$$
\mult_{P_i}( \Qd.\yczt) = \mbox{i}(S.\yczt.H;P_i) =   \mbox{i}( (S \cap H). (\yczt \cap H);P_i) = 2 \, ,
$$  
where the last equality results from the fact that the
tangent cone $\mbox{C}_{P_i}\yczt$ consists of the  planes $\mbox{T}_{P_i} \sj, \mbox{T}_{P_i} \sd$,
so the planar curves  $\Qd \cap H$, $\yczt \cap H$ have no common tangent lines at $P_i$.
The curves  $C'_4, C''_4$ are components of the cycle
$\Qd.\yczt$ and pass through the point $P_i$, 
so they are smooth in the cusp.

\vskip3mm
\noindent
{\bf Example~\ref{examples}.2} {\rm (Configuration (II))} We put
$$
L' := x_0, \, L'' := x_1, \,  F' := x_2, \,  F'' :=  6 (x_1 + x_2) - 11 x_0,
$$
and define  a quartic $\yczt$ by the equation   (\ref{determinantalequation}) with 
$\Qd := x_3^2 - x_2^2.$ \\
Then the curve $C_3$ is a sum of the lines $l_1, l_2, l_3$, where 
$l_j : x_0 = j x_2 , \,  x_1 = j^2 x_2$, and $\Qd$ meets $C_3$ in the points
$( j : j^2 : 1 : \pm 1 ), \mbox{ where } j=1,2,3.$ \\
If we pass to the equation (\ref{minus}), then
we obtain the residual quadric  
$$
R \, : \, x_3^2 - x_2^2 - x_0 \cdot x_1. 
$$
One can check that $\sj, \sd, \Qd$ meet transversally in the six points, none of which 
belongs to $R$, so those points are cusps on $\yczt$. 
A Gr\"obner basis computation shows that the quartic  $\yczt$ has no extra singularities.
Finally, one can imitate the proof in Example~\ref{examples}.1 to show that the set $\sing( \yczt )$ is $3$-divisible. We leave the 
details to the reader.

\vskip2mm
For the sake of completeness we explain below the way we applied Gr\"obner bases to 
count singularities of quartics.

\vskip2mm
\noindent
{\bf Remark~\ref{examples}.1} 
In order to show that a polynomial $g \in \kk[x_0, x_1, x_2, x_3]$
vanishes on
 all singularities of a surface, 
we  need the notion of  the 
 remainder on division 
of a polynomial  $g$
 by a  Gr\"obner basis  ${\cal B}$   of an ideal 
${\cal I} \subset \kk[x_0, x_1, x_2, x_3]$. 
Its definition can be found in \cite[II.$\S$6]{clos}.
The only fact we need is that if the remainder  
vanishes, then
the polynomial $g$ belongs to the ideal ${\cal I}$ (see \cite[II.$\S$6.Cor.~2]{clos}).
We use this fact 
in the following way: \\
Let  ${\cal B}$ be 
a  Gr\"obner basis
of the jacobian ideal of the quartic  $Y_4$.
If we find such an integer $p$ that 
the remainder on division of the polynomial $g^p$ by the basis ${\cal B}$  vanishes, 
then all singularities of the surface $Y_4$ lie on the hypersurface given by 
$g$. 
The former can checked with the following Maple~6.0 commands,
where  Yi is defined as  the partial derivative $\frac{\partial Y_4}{\partial x_i}$:
\begin{verbatim} 
with(Groebner):
B := gbasis([Y0, Y1, Y2, Y3], tdeg( x0, x1, x2, x3 )):
normalf((g^p), B,  tdeg( x0, x1, x2, x3 )); 
\end{verbatim}  
If the output is zero, then 
$g^p$ belongs to the jacobian ideal of $Y_4$. 

\vskip2mm
\noindent
{\bf Acknowledgements:} I would like to express my sincere gratitude to Prof. W. Barth 
for his constant help and numerous valuable discussions. 


\noindent
{\small Current address: \\
Mathematisches Institut, FAU Erlangen-N\"urnberg, Bismarckstrasse 1 1/2, \\
D-91054 Erlangen, Germany \\
Permanent address: \\
 Institute of Mathematics, Jagiellonian University, ul.~Reymonya~4, \\
30-059 Krak\'ow, Poland \\
{\sl E-mail address:} rams@mi.uni-erlangen.de, rams@im.uj.edu.pl}

\end{document}